\documentclass[a4paper,11pt,oneside]{article}

\usepackage{amsmath,amsthm,amsfonts,amssymb,amscd}
\usepackage{mathtools}
\usepackage{pgf,tikz}
\usepackage{float}
\usepackage{enumerate}
\usepackage{booktabs}
\usepackage{pgfplots}
\usepackage[noend]{algpseudocode}
\usepackage{graphicx,bm,xcolor}
\usepackage{algorithm}
\usepackage{algpseudocode}
\usepackage{hyperref}
\usepackage{caption}
\usepackage[T1]{fontenc}
\usepackage[polish,german,english]{babel}
\usepackage{verbatim} 

\usepackage{url}
\usepackage{authblk}

\usepackage[top=2cm, bottom=2cm, left=2cm, right=2cm]{geometry}

\theoremstyle{plain} 
\newtheorem{theorem}{Theorem}[section]

\newtheorem{remark}[theorem]{Remark}
\newtheorem{definition}[theorem]{Definition}

\theoremstyle{definition} %

\theoremstyle{remark} %

\newcommand{\R}{\mathbb R} 
\newcommand{\C}{\mathbb C}
\newcommand{\dd}{\mathsf{d}}
\usepackage{yfonts}
\usepackage{dsfont} 
\usepackage[breakable]{tcolorbox}
\usepackage{upquote}
\usepackage{fancyvrb}
\definecolor{urlcolor}{rgb}{0,.145,.698}
\definecolor{linkcolor}{rgb}{.71,0.21,0.01}
\definecolor{citecolor}{rgb}{.12,.54,.11}
\definecolor{ansi-black}{HTML}{3E424D}
\definecolor{ansi-black-intense}{HTML}{282C36}
\definecolor{ansi-red}{HTML}{E75C58}
\definecolor{ansi-red-intense}{HTML}{B22B31}
\definecolor{ansi-green}{HTML}{00A250}
\definecolor{ansi-green-intense}{HTML}{007427}
\definecolor{ansi-yellow}{HTML}{DDB62B}
\definecolor{ansi-yellow-intense}{HTML}{B27D12}
\definecolor{ansi-blue}{HTML}{208FFB}
\definecolor{ansi-blue-intense}{HTML}{0065CA}
\definecolor{ansi-magenta}{HTML}{D160C4}
\definecolor{ansi-magenta-intense}{HTML}{A03196}
\definecolor{ansi-cyan}{HTML}{60C6C8}
\definecolor{ansi-cyan-intense}{HTML}{258F8F}
\definecolor{ansi-white}{HTML}{C5C1B4}\definecolor{ansi-white-intense}{HTML}{A1A6B2}
\definecolor{ansi-default-inverse-fg}{HTML}{FFFFFF}
\definecolor{ansi-default-inverse-bg}{HTML}{000000}
\definecolor{incolor}{HTML}{303F9F}
\definecolor{outcolor}{HTML}{D84315}
\definecolor{cellborder}{HTML}{CFCFCF}
\definecolor{cellbackground}{HTML}{F7F7F7}

\makeatletter
\def\PY@reset{\let\PY@it=\relax \let\PY@bf=\relax%
    \let\PY@ul=\relax \let\PY@tc=\relax%
    \let\PY@bc=\relax \let\PY@ff=\relax}
\def\PY@tok#1{\csname PY@tok@#1\endcsname}
\def\PY@toks#1+{\ifx\relax#1\empty\else%
    \PY@tok{#1}\expandafter\PY@toks\fi}
\def\PY@do#1{\PY@bc{\PY@tc{\PY@ul{%
    \PY@it{\PY@bf{\PY@ff{#1}}}}}}}
\def\PY#1#2{\PY@reset\PY@toks#1+\relax+\PY@do{#2}}

\expandafter\def\csname PY@tok@w\endcsname{\def\PY@tc##1{\textcolor[rgb]{0.73,0.73,0.73}{##1}}}
\expandafter\def\csname PY@tok@c\endcsname{\let\PY@it=\textit\def\PY@tc##1{\textcolor[rgb]{0.25,0.50,0.50}{##1}}}
\expandafter\def\csname PY@tok@cp\endcsname{\def\PY@tc##1{\textcolor[rgb]{0.74,0.48,0.00}{##1}}}
\expandafter\def\csname PY@tok@k\endcsname{\let\PY@bf=\textbf\def\PY@tc##1{\textcolor[rgb]{0.00,0.50,0.00}{##1}}}
\expandafter\def\csname PY@tok@kp\endcsname{\def\PY@tc##1{\textcolor[rgb]{0.00,0.50,0.00}{##1}}}
\expandafter\def\csname PY@tok@kt\endcsname{\def\PY@tc##1{\textcolor[rgb]{0.69,0.00,0.25}{##1}}}
\expandafter\def\csname PY@tok@o\endcsname{\def\PY@tc##1{\textcolor[rgb]{0.40,0.40,0.40}{##1}}}
\expandafter\def\csname PY@tok@ow\endcsname{\let\PY@bf=\textbf\def\PY@tc##1{\textcolor[rgb]{0.67,0.13,1.00}{##1}}}
\expandafter\def\csname PY@tok@nb\endcsname{\def\PY@tc##1{\textcolor[rgb]{0.00,0.50,0.00}{##1}}}
\expandafter\def\csname PY@tok@nf\endcsname{\def\PY@tc##1{\textcolor[rgb]{0.00,0.00,1.00}{##1}}}
\expandafter\def\csname PY@tok@nc\endcsname{\let\PY@bf=\textbf\def\PY@tc##1{\textcolor[rgb]{0.00,0.00,1.00}{##1}}}
\expandafter\def\csname PY@tok@nn\endcsname{\let\PY@bf=\textbf\def\PY@tc##1{\textcolor[rgb]{0.00,0.00,1.00}{##1}}}
\expandafter\def\csname PY@tok@ne\endcsname{\let\PY@bf=\textbf\def\PY@tc##1{\textcolor[rgb]{0.82,0.25,0.23}{##1}}}
\expandafter\def\csname PY@tok@nv\endcsname{\def\PY@tc##1{\textcolor[rgb]{0.10,0.09,0.49}{##1}}}
\expandafter\def\csname PY@tok@no\endcsname{\def\PY@tc##1{\textcolor[rgb]{0.53,0.00,0.00}{##1}}}
\expandafter\def\csname PY@tok@nl\endcsname{\def\PY@tc##1{\textcolor[rgb]{0.63,0.63,0.00}{##1}}}
\expandafter\def\csname PY@tok@ni\endcsname{\let\PY@bf=\textbf\def\PY@tc##1{\textcolor[rgb]{0.60,0.60,0.60}{##1}}}
\expandafter\def\csname PY@tok@na\endcsname{\def\PY@tc##1{\textcolor[rgb]{0.49,0.56,0.16}{##1}}}
\expandafter\def\csname PY@tok@nt\endcsname{\let\PY@bf=\textbf\def\PY@tc##1{\textcolor[rgb]{0.00,0.50,0.00}{##1}}}
\expandafter\def\csname PY@tok@nd\endcsname{\def\PY@tc##1{\textcolor[rgb]{0.67,0.13,1.00}{##1}}}
\expandafter\def\csname PY@tok@s\endcsname{\def\PY@tc##1{\textcolor[rgb]{0.73,0.13,0.13}{##1}}}
\expandafter\def\csname PY@tok@sd\endcsname{\let\PY@it=\textit\def\PY@tc##1{\textcolor[rgb]{0.73,0.13,0.13}{##1}}}
\expandafter\def\csname PY@tok@si\endcsname{\let\PY@bf=\textbf\def\PY@tc##1{\textcolor[rgb]{0.73,0.40,0.53}{##1}}}
\expandafter\def\csname PY@tok@se\endcsname{\let\PY@bf=\textbf\def\PY@tc##1{\textcolor[rgb]{0.73,0.40,0.13}{##1}}}
\expandafter\def\csname PY@tok@sr\endcsname{\def\PY@tc##1{\textcolor[rgb]{0.73,0.40,0.53}{##1}}}
\expandafter\def\csname PY@tok@ss\endcsname{\def\PY@tc##1{\textcolor[rgb]{0.10,0.09,0.49}{##1}}}
\expandafter\def\csname PY@tok@sx\endcsname{\def\PY@tc##1{\textcolor[rgb]{0.00,0.50,0.00}{##1}}}
\expandafter\def\csname PY@tok@m\endcsname{\def\PY@tc##1{\textcolor[rgb]{0.40,0.40,0.40}{##1}}}
\expandafter\def\csname PY@tok@gh\endcsname{\let\PY@bf=\textbf\def\PY@tc##1{\textcolor[rgb]{0.00,0.00,0.50}{##1}}}
\expandafter\def\csname PY@tok@gu\endcsname{\let\PY@bf=\textbf\def\PY@tc##1{\textcolor[rgb]{0.50,0.00,0.50}{##1}}}
\expandafter\def\csname PY@tok@gd\endcsname{\def\PY@tc##1{\textcolor[rgb]{0.63,0.00,0.00}{##1}}}
\expandafter\def\csname PY@tok@gi\endcsname{\def\PY@tc##1{\textcolor[rgb]{0.00,0.63,0.00}{##1}}}
\expandafter\def\csname PY@tok@gr\endcsname{\def\PY@tc##1{\textcolor[rgb]{1.00,0.00,0.00}{##1}}}
\expandafter\def\csname PY@tok@ge\endcsname{\let\PY@it=\textit}
\expandafter\def\csname PY@tok@gs\endcsname{\let\PY@bf=\textbf}
\expandafter\def\csname PY@tok@gp\endcsname{\let\PY@bf=\textbf\def\PY@tc##1{\textcolor[rgb]{0.00,0.00,0.50}{##1}}}
\expandafter\def\csname PY@tok@go\endcsname{\def\PY@tc##1{\textcolor[rgb]{0.53,0.53,0.53}{##1}}}
\expandafter\def\csname PY@tok@gt\endcsname{\def\PY@tc##1{\textcolor[rgb]{0.00,0.27,0.87}{##1}}}
\expandafter\def\csname PY@tok@err\endcsname{\def\PY@bc##1{\setlength{\fboxsep}{0pt}\fcolorbox[rgb]{1.00,0.00,0.00}{1,1,1}{\strut ##1}}}
\expandafter\def\csname PY@tok@kc\endcsname{\let\PY@bf=\textbf\def\PY@tc##1{\textcolor[rgb]{0.00,0.50,0.00}{##1}}}
\expandafter\def\csname PY@tok@kd\endcsname{\let\PY@bf=\textbf\def\PY@tc##1{\textcolor[rgb]{0.00,0.50,0.00}{##1}}}
\expandafter\def\csname PY@tok@kn\endcsname{\let\PY@bf=\textbf\def\PY@tc##1{\textcolor[rgb]{0.00,0.50,0.00}{##1}}}
\expandafter\def\csname PY@tok@kr\endcsname{\let\PY@bf=\textbf\def\PY@tc##1{\textcolor[rgb]{0.00,0.50,0.00}{##1}}}
\expandafter\def\csname PY@tok@bp\endcsname{\def\PY@tc##1{\textcolor[rgb]{0.00,0.50,0.00}{##1}}}
\expandafter\def\csname PY@tok@fm\endcsname{\def\PY@tc##1{\textcolor[rgb]{0.00,0.00,1.00}{##1}}}
\expandafter\def\csname PY@tok@vc\endcsname{\def\PY@tc##1{\textcolor[rgb]{0.10,0.09,0.49}{##1}}}
\expandafter\def\csname PY@tok@vg\endcsname{\def\PY@tc##1{\textcolor[rgb]{0.10,0.09,0.49}{##1}}}
\expandafter\def\csname PY@tok@vi\endcsname{\def\PY@tc##1{\textcolor[rgb]{0.10,0.09,0.49}{##1}}}
\expandafter\def\csname PY@tok@vm\endcsname{\def\PY@tc##1{\textcolor[rgb]{0.10,0.09,0.49}{##1}}}
\expandafter\def\csname PY@tok@sa\endcsname{\def\PY@tc##1{\textcolor[rgb]{0.73,0.13,0.13}{##1}}}
\expandafter\def\csname PY@tok@sb\endcsname{\def\PY@tc##1{\textcolor[rgb]{0.73,0.13,0.13}{##1}}}
\expandafter\def\csname PY@tok@sc\endcsname{\def\PY@tc##1{\textcolor[rgb]{0.73,0.13,0.13}{##1}}}
\expandafter\def\csname PY@tok@dl\endcsname{\def\PY@tc##1{\textcolor[rgb]{0.73,0.13,0.13}{##1}}}
\expandafter\def\csname PY@tok@s2\endcsname{\def\PY@tc##1{\textcolor[rgb]{0.73,0.13,0.13}{##1}}}
\expandafter\def\csname PY@tok@sh\endcsname{\def\PY@tc##1{\textcolor[rgb]{0.73,0.13,0.13}{##1}}}
\expandafter\def\csname PY@tok@s1\endcsname{\def\PY@tc##1{\textcolor[rgb]{0.73,0.13,0.13}{##1}}}
\expandafter\def\csname PY@tok@mb\endcsname{\def\PY@tc##1{\textcolor[rgb]{0.40,0.40,0.40}{##1}}}
\expandafter\def\csname PY@tok@mf\endcsname{\def\PY@tc##1{\textcolor[rgb]{0.40,0.40,0.40}{##1}}}
\expandafter\def\csname PY@tok@mh\endcsname{\def\PY@tc##1{\textcolor[rgb]{0.40,0.40,0.40}{##1}}}
\expandafter\def\csname PY@tok@mi\endcsname{\def\PY@tc##1{\textcolor[rgb]{0.40,0.40,0.40}{##1}}}
\expandafter\def\csname PY@tok@il\endcsname{\def\PY@tc##1{\textcolor[rgb]{0.40,0.40,0.40}{##1}}}
\expandafter\def\csname PY@tok@mo\endcsname{\def\PY@tc##1{\textcolor[rgb]{0.40,0.40,0.40}{##1}}}
\expandafter\def\csname PY@tok@ch\endcsname{\let\PY@it=\textit\def\PY@tc##1{\textcolor[rgb]{0.25,0.50,0.50}{##1}}}
\expandafter\def\csname PY@tok@cm\endcsname{\let\PY@it=\textit\def\PY@tc##1{\textcolor[rgb]{0.25,0.50,0.50}{##1}}}
\expandafter\def\csname PY@tok@cpf\endcsname{\let\PY@it=\textit\def\PY@tc##1{\textcolor[rgb]{0.25,0.50,0.50}{##1}}}
\expandafter\def\csname PY@tok@c1\endcsname{\let\PY@it=\textit\def\PY@tc##1{\textcolor[rgb]{0.25,0.50,0.50}{##1}}}
\expandafter\def\csname PY@tok@cs\endcsname{\let\PY@it=\textit\def\PY@tc##1{\textcolor[rgb]{0.25,0.50,0.50}{##1}}}

\def\PYZus{\char`\_}

\def\PYZlt{\char`\<}
\def\PYZgt{\char`\>}
\def\PYZsh{\char`\#}
\def\PYZpc{\char`\%}

\def\PYZhy{\char`\-}

\def\PYZdq{\char`\"}

\DefineVerbatimEnvironment{Highlighting}{Verbatim}{commandchars=\\\{\}}

\begin{document}

\title{Domain Decomposition with Neural Network Interface Approximations for time-harmonic Maxwell's equations with different wave numbers}

\author[1]{T. Knoke}
\author[1,3]{S. Kinnewig}
\author[1,3]{S. Beuchler}
\author[2,3]{A. Demircan}
\author[2,3]{U. Morgner}
\author[1,3]{T. Wick}

\affil[1]{Leibniz University Hannover, Institute of Applied Mathematics, Germany}
\affil[2]{Leibniz University Hannover, Institute of Quantum Optics, Germany }
\affil[3]{Leibniz University Hannover, Cluster of Excellence PhoenixD (Photonics, Optics, and
        Engineering - Innovation Across Disciplines), Germany}

\date{}

\maketitle
	
\begin{abstract}
In this work, we consider the time-harmonic Maxwell’s equations and their 
  numerical solution with a domain decomposition method. As an innovative feature, 
  we propose a feedforward neural network-enhanced approximation of the interface conditions 
  between the subdomains. The advantage is that the interface condition 
can be updated without recomputing the Maxwell system at each step.
The main part consists of a detailed description of the construction 
of the neural network for domain decomposition and the training process.
To substantiate this proof of concept, we
  investigate a few subdomains in some numerical experiments with low 
  frequencies. Therein the new approach is compared to a classical domain decomposition method.
Moreover, we highlight 
  current challenges of training and testing with different wave numbers 
  and we provide information on the behaviour of the 
  neural-network, such as convergence of the loss function, 
  and different activation functions.
\end{abstract}

\begin{keywords}
{Time-Harmonic Maxwell's Equations, Machine Learning, Feedforward Neural Network, Domain Decomposition Method.}
\end{keywords}	
	
\section{Introduction}
\label{sec_intro}
The Maxwell’s equations for describing electro-magnetic phenomena are of great interest in current research fields, such as optics. One present 
example of employing Maxwell's equations 
can be found in the Cluster of Excellence PhoenixD (Photonics Optics Engineering Innovation Across Disciplines)\footnote{\url{https://www.phoenixd.uni-hannover.de/en/}} at the Leibniz University Hannover,
in which modern methods for optics simulations are being developed. 
Therein, one focus
is on the efficient and accurate calculation of light distribution in an optical material to design optical devices on the micro- and nanoscale \cite{Shi,Melchert}. 
In comparison to other partial differential equations, such as in solid 
mechanics or fluid flow, the Maxwell's equations have some peculiarities
such as the curl operator, which has in two-dimensional problems, a one-dimensional image, but 
in three-dimensional problems, it has a three-dimensional image. Moreover, 
the requirements for the discretization and definiteness of the final 
linear system are specific. In more detail, 
in numerical mathematics, Maxwell’s equations are of interest because 
of their specific mathematical structures \cite{bk:monk:03,Dem06,LaPauRe19,AlonsoRodriguezBertolazziValli+2019+1+44}, 
requirements for finite elements 
\cite{bk:monk:03,LaPauRe19,art:nedelec:80,art:nedelec:86,art:kinnewig:23,CARSTENSEN2016494,NicaiseTomezyk+2019+285+340}, 
their numerical solution \cite{Hip98,HENNEKING202185,faustmann2022h,art:dol09} 
as well as postprocessing such as a posteriori error control and adaptivity \cite{Schoe08,Bue12}.
As their numerical solution is challenging due to their ill-posed nature, e.g., \cite{art:beu:21}, 
one must apply suitable techniques. 
The most prominent approach in the literature is based 
on domain decomposition (DD) techniques \cite{bk:toselli, art:dol09}. 
The geometric multigrid solver developed by Hiptmair \cite{Hip98} can only be 
applied to the problem in the time domain 
(i.e., the well-posed problem).

In this work, we concentrate on the numerical solution 
using a domain decomposition 
method. Specifically, our starting point is the method developed in \cite{art:beu:21}, 
based on ideas from \cite{art:bour:15},
and which was realized in the modern open-source finite element library 
deal.II \cite{deal2020,software:dealII94}.
The domain decomposition method's crucial point is the interface operator derivation \cite{art:bour:15}.
Our main objective in the current work is to 
design a proof of concept to approximate 
the interface operator with the 
help of a feedforward neural network (NN) \cite{Bish06,HiHi19,bk:KiKoRoWi:22}.
We carefully derive the governing algorithms and focus on a two-domain problem 
to study our new approach's mechanism and performance.
Implementation-wise, the previously mentioned deal.II library (in C++) is 
coupled to the PyTorch (in python) \cite{software:pytorch} library, which is one of the 
standard packages for neural network computations.
Our main aim is to showcase that our approach is 
feasible and can be a point of departure for 
further future extensions. We notice that the current work
is an extension of the conference proceedings paper \cite{KnoKiBeuWi22_dd27}
with more mathematical and algorithmic details, and different numerical 
tests, specifically the studies on different wave numbers and comparison 
of two NN activation functions.

The outline of this work is as follows: In Section \ref{sec_equations}, we introduce the 
time-harmonic Maxwell's equations and our notation.
Next, in Section \ref{sec_numerics}, domain decomposition and neural network 
approximations are introduced.
Afterwards, we address in detail the training process in Section \ref{sec_training}.
In Section \ref{sec_tests}, some numerical tests demonstrate our proof of concept. 
Our work is summarized in Section \ref{sec_conclusions}.

\section{Equations}
\label{sec_equations}
For the sake of simplicity, we only consider the two-dimensional time-harmonic Maxwell's 
equations. In the following, we will introduce these equations in detail.

\subsection{Fundamental operators}
To comprehensively describe the problem, we introduce the basic operators needed 
to describe two-dimensional electro-magnetic problems. Therefore, let us assume 
a scalar function $\phi: \R \rightarrow \R$ and $\vec v \in \R^2$ to be a 
two-dimensional vector. Then the gradient of $\phi$ is given by
$
  \nabla \phi = \left( \frac{\partial \phi}{\partial x_1},~\frac{\partial \phi}{\partial x_2} \right),
$
and the divergence of $v$ is given by
$
 \operatorname{div}(v) \coloneqq \nabla \cdot v \coloneqq \sum_{i = 1}^{2} \frac{ \partial v_i }{ \partial x_i }.
$
Next,
$a \cdot b = (a_1, a_2)^T \cdot (b_1, b_2)^T = a_1 b_1 + a_2 b_2$ denotes the
scalar product.
We can furthermore write down the description of the two-dimensional curl operator
\begin{equation}
\operatorname{curl}(\vec v)=  \frac{\partial v_2}{\partial x_1} - \frac{\partial v_1}{\partial x_2},
\end{equation}
and the curl operator applied to a scalar function
\begin{equation}
  \underline{\operatorname{curl}}(\phi) = \left( \begin{array}{c} \frac{ \partial \phi }{ \partial x_2 } \\ - \frac{ \partial \phi }{ \partial x_1 } \end{array} \right).
\end{equation}

\subsection{Time-harmonic Maxwell's equations}
Let $\Omega \subset \R^2$ be a bounded domain with sufficiently smooth boundary $\Gamma$. 
The latter is partitioned into $\Gamma=\Gamma^{\infty} \cup \Gamma^{\text{inc}}$. 
The main governing function space is defined as 
\[
H(\mathrm{curl},\Omega):=\{v \in \mathcal{L}^2(\Omega) \ | \ \mathrm{curl}(v) \in \mathcal{L}^2(\Omega)\}, 
\]
where $\mathcal{L}^2(\Omega)$ is the well-known space of square-integrable 
functions in the Lebesgue sense. In order to define boundary conditions, 
we introduce the traces 
\begin{align*}
&\gamma^t: H(\mathrm{curl},\Omega) \to H_{\times}^{-1/2}(\mathrm{div},\Gamma),\\ 
&\gamma^T: H(\mathrm{curl},\Omega) \to H_{\times}^{-1/2}(\mathrm{curl},\Gamma),
\end{align*}
which are defined by
\begin{align*}
  \gamma^t\left(\phi\right)= \left( \begin{array}{c} ~ \phi \; n_2 \\ - \phi \; n_1 \end{array} \right) \quad \text{and} \quad
\gamma^T\left(v\right)
= v - (n \cdot v) \cdot n,
\end{align*}
where $n \in \R^2$ is the normal vector of 
$\Omega$, $H_{\times}^{-1/2}(\mathrm{div},\Gamma):=\{v \in H^{-1/2}(\Gamma) \ | \ v \cdot n=0, \ \mathrm{div}_{\Gamma} v \in H^{-1/2}(\Gamma)\}$ 
is the space of well-defined surface divergence fields and 
$H(\mathrm{curl},\Gamma):=\{v \in H^{-1/2}(\Gamma) \ | \ v \cdot n=0, \ \mathrm{curl}_{\Gamma} \left( v \right) \in H^{-1/2}(\Gamma)\}$ 
is the space of well-defined surface curls, see \cite[Chapter 3.4]{bk:monk:03}.
In the following, we first state the strong form of the system.
The time-harmonic Maxwell's equations are then defined as follows: 
Find the electric field $E:\Omega\to\C^2$ such that
\begin{align}
  \left\{
  \begin{array}{l l l}
    \underline{\operatorname{curl}}\left( \mu^{-1} \operatorname{curl}\left(\vec E \right) \right) - \varepsilon \omega^2 \vec E            & = \vec 0            & \text{ in } \Omega          \\
    \mu^{-1} \gamma^t \left( \operatorname{curl}\left( \vec E \right) \right) - i \kappa \omega \gamma^T \left( \vec E \right)  & = \vec 0            & \text{ on } \Gamma^{\infty} \\
    \gamma^T\left(\vec E\right)                                                                                                 & = \vec E^\text{inc} & \text{ on } \Gamma^{\text{inc}},
    \label{maxwellsystem}
  \end{array}
  \right.
\end{align}
where $\vec E^\text{inc}:\mathbb{R}^2\rightarrow\mathbb{C}^2$ is some given incident electric field, 
$\mu \in \mathbb{R}^+$ is the relative magnetic permeability, $\kappa = \sqrt{\varepsilon}$, $\varepsilon \in \mathbb{C}$ 
relative permittivity, $\omega = \frac{2 \pi}{\lambda}$ is the wave number and $\lambda \in \mathbb{R}^{+}$ is the wave length
and $i$ denotes the imaginary number. 
System (\ref{maxwellsystem}), as well as its weak form, 
is called time-harmonic, because the time dependence can be expressed by 
$e^{i \omega \tau}$, where $\tau \ge 0$ denotes the time. 

\subsection{Weak formulation}
In this subsection, we derive the weak form. This is the 
starting point for a finite element method (FEM) discretization. 
For the derivation, we first begin by rewriting 
the curl product with the help of integration by parts:
\begin{align}
  \label{eq:integration_by_parts}
  \int_\Omega \underline{\operatorname{curl}} \left( \phi \right) \cdot \vec u ~ \dd x 
    = \int_\Omega \phi \operatorname{curl} \left( u \right) ~ \dd x
    + \int_{\partial \Omega} \gamma^t(\phi) \cdot u ~ \dd s,
\end{align}
see for instance \cite{GirRav,bk:monk:03}.
We want to derive the weak formulation from the strong formulation 
\eqref{maxwellsystem}
in the following:
\begin{align}
  & \int_\Omega \underline{\operatorname{curl}}\left(\mu^{-1} \operatorname{curl} \left( \vec E \right) \right) \cdot \vec \varphi ~ \dd x
    - \varepsilon \omega^2 \int_\Omega \vec E \cdot \vec \varphi ~ \dd x
    = \vec 0, 
  \nonumber \\
  \overset{\text{\eqref{eq:integration_by_parts}}}{\Rightarrow}
  & \int_\Omega \mu^{-1} \operatorname{curl} \left( \vec E \right) \operatorname{curl} \left( \vec \varphi \right) ~ \dd x
    - \varepsilon \omega^2 \int_\Omega \vec E \cdot \vec \varphi ~ \dd x
    + \int_{\partial \Omega} \mu^{-1} \gamma^t \left( \operatorname{curl} \left( \vec E \right) \right) \cdot \vec \varphi ~ \dd s
    = \vec 0.
  \label{eq:thm_derivation}
\end{align}
By applying the definition of the boundaries $\Gamma^\infty$ and $\Gamma^\text{inc}$ from 
equation 
\eqref{maxwellsystem}
to equation \eqref{eq:thm_derivation}, we obtain the weak 
formulation of the time-harmonic Maxwell's equations.
Find $\vec E \in H(curl,\Omega)$ such that for all 
$\vec \varphi \in H(curl,\Omega)$
\begin{align}
  \int_\Omega \left( \mu^{-1} \operatorname{curl} \left( \vec E \right) \operatorname{curl} \left( \vec \varphi \right) 
    - \varepsilon \omega^2 \vec E \cdot \vec \varphi \right) ~ \dd x 
    + & i \kappa \omega \int_{\Gamma^\infty} \gamma^T \left( \vec E \right) \cdot \gamma^T \left( \vec \varphi \right) ~ \dd s 
    \nonumber \\ 
    & = \int_{\Gamma^\text{inc}} \gamma^T \left( \vec E^\text{inc} \right) \cdot \gamma^T \left( \vec \varphi \right) ~ \dd s.
  \label{eq:thm_weak}
\end{align}

\subsection{Two-dimensional N\'ed\'elec elements}
For the implementation with the help of a Galerkin finite element method (FEM), 
we need the discrete weak form. Based on the De-Rham cohomology, we must choose 
our basis functions out of the N\'ed\'elec space $V_h$. Therefore, we want to 
introduce the definition of the space $V_h$ in the following, based on the 
formalism introduced by Zaglmayr \cite[Chapter 5.2]{bk:zaglmayr}. 

As a suitable polynomial basis, we introduce the integrated Legendre polynomials.
Let $x \in [-1,1]$.
The following recursive formula defines the integrated Legendre polynomials:
\begin{equation}
  \begin{array}{rl}
    L_1(x) & = x, \\[1mm]
    L_2(x) & = \frac{1}{2}\left( x^2 - 1 \right), \\[1mm]
    (n + 1) L_{n+1}(x) & = (2n - 1)x L_n(x) - (n - 2) L_{n-1}(x), \quad\text{ for } n \geq 2.
  \end{array}
\end{equation}
Let us choose the quadrilateral reference element as $Q = [0, 1] \times [0, 1]$.
\begin{figure}[H]
\begin{minipage}{0.25\textwidth}
  \begin{tikzpicture}[]
    \coordinate (V0) at (0.0, 0.0);
    \coordinate (V1) at (3.0, 0.0);
    \coordinate (V2) at (0.0, 3.0);
    \coordinate (V3) at (3.0, 3.0);
  
    \foreach \n in {V0, V1, V2, V3}
      \node at (\n)[circle, fill, inner sep=1.5pt]{};
  
    \node[anchor=north west] at (V0) {$V_0$};
    \node[below] at (V1) {$V_1$};
    \node[anchor=south east] at (V2) {$V_2$};
    \node[above] at (V3) {$V_3$};
  
    \node[anchor=east, gray] at (-0.2, 1.5) {$\mathcal{E}_0$};
    \node[anchor=west, gray] at (3.2, 1.5) {$\mathcal{E}_1$};
    \node[anchor=north, gray] at (1.5, -0.2) {$\mathcal{E}_2$};
    \node[anchor=south, gray] at (1.5, 3.2) {$\mathcal{E}_3$};
  
  
    \coordinate (Vx) at (3.7, 0.0);
    \coordinate (Vy) at (0.0, 3.7);
    
    \draw[->, >=stealth] (-0.5, 0) -> (Vx);
    \draw[->, >=stealth] (0, -0.5) -> (Vy);
    
    \node[right] at (Vx) {x};
    \node[above] at (Vy) {y};
  
    \draw[thick]      (V0) -- (V1);
    \draw[thick]      (V0) -- (V2);
    \draw[thick]      (V3) -- (V2);
    \draw[thick]      (V3) -- (V1);
  
    \draw[->, >=stealth, thin, gray] (-0.2, 0.5) -- (-0.2, 2.5);
    \draw[->, >=stealth, thin, gray] (3.2, 0.5) -- (3.2, 2.5);
    \draw[->, >=stealth, thin, gray] (0.5, -0.2) -- (2.5, -0.2);
    \draw[->, >=stealth, thin, gray] (0.5, 3.2) -- (2.5, 3.2);
  \end{tikzpicture}
\end{minipage}
\hfill
\begin{minipage}{0.65\textwidth}
    \[
    \begin{array}{l l}
      \lambda_0 = (1-x) (1-y), & \sigma_0 = (1-x) + (1-y), \\ 
      \lambda_1 =    x  (1-y), & \sigma_1 =    x  + (1-y), \\
      \lambda_2 = (1-x)    y , & \sigma_2 = (1-x) +    y , \\
      \lambda_3 =    x     y , & \sigma_3 =    x  +    y  \\
    \end{array}
  \]
\end{minipage} 
  \caption{
    \label{fig:reference_cell}
    Left: Vertex and edge ordering on the reference cell, right: parametrisation 
    of the reference cell.
  }
\end{figure}

We continue by defining the set of all edges 
$\textgoth{E} = \left\{ \mathcal{E}_m \right\}_{0 \leq m < 4}$ 
with local edge-ordering $\mathcal{E}_m = \{ V_i, V_j \}$, where $(i,j)\in\{(0,2),(1,3),(0,1),(2,3) \}$, see Figure \ref{fig:reference_cell}. We  denote the 
cell itself with local vertex-ordering $C = \{V_0, V_1, V_2, V_3\}$. The polynomial 
order is given by $\vec p = \left( \left\{ p_\mathcal{E} \right\}_{\mathcal{E} \in \textgoth{E}}, p_C \right)$.
\begin{figure}[H]
\begin{center}
  \begin{tabular} {|p{3cm}|p{10cm}|} \hline
    \multicolumn{2}{|c|}{\textbf{ $H(curl)$ conforming basis function}} \\ \hline \hline
    \multicolumn{2}{|l|}{\textbf{Vertex-based shape functions}}                        \\ \hline
    \multicolumn{2}{|l|}{There are no DoFs on the vertices.}                           \\ \hline \hline
    \multicolumn{2}{|l|}{\textbf{Edge-based shape functions}}       \\
    \multicolumn{2}{|l|}{for $0 \leq i < p_\mathcal{E}$, $\mathcal{E} \in \textgoth{E}$, where $\lambda_\alpha$ and $\sigma_\alpha,~\alpha \in \{0,1,2,3\}$ are defined in Figure \ref{fig:reference_cell}} \\ \hline
    Lowest order &
    $
      \varphi_{\mathcal{E}_m}^{\mathcal{N}_0} = 
      \frac{1}{2} \nabla \left( \sigma_{e_2} - \sigma_{e_1} \right) \left( \lambda_{e_1} + \lambda_{e_2} \right)
    $ \\[1mm] \hline
    Higher-order &
    $
      \varphi_{i}^{\mathcal{E}_m} = 
      \nabla \left( L_{i+2} \left( \sigma_{e_2} - \sigma_{e_1} \right) \left( \lambda_{e_1} + \lambda_{e_2} \right) \right)
    $ \\[1mm] \hline \hline
    
    \multicolumn{2}{|l|}{\textbf{Cell-based functions}}       \\
    \multicolumn{2}{|l|}{$0 \leq i,j < p_C$}                  \\[1mm] \hline
    Type 1: &
    $
      \varphi_{(i,j)}^{C,1} = \nabla ( L_{i+2}(\xi_F) L_{j+2}(\eta_F) )
    $ \\[1mm] \hline
    Type 2: &
    $
      \varphi_{(i,j)}^{C,2} = \widetilde \nabla ( L_{i+2}(\xi_F) L_{j+2}(\eta_F) )
    $ \\[1mm]
    & ~~where $\widetilde \nabla(a~b) := (a' ~ b - a ~ b')$ \\ \hline
    Type 3: &
    $
      \varphi_{(0,j)}^{C,3} = L_{i+2}(2y-1) \vec e_x 
    $ \\[1mm]
    & $
      \varphi_{(i,0)}^{C,3} = L_{i+2}(2x-1) \vec e_y 
    $ \\[1mm] \hline
  \end{tabular}
  \caption{
    \label{fig:nedelec_space}
    The definition of the ${H}(\operatorname{curl})$ basis-functions on the reference element with barycentric coordinates $\lambda_\alpha$ and $\sigma_\alpha,~\alpha \in \{0,1,2,3\}$. 
  }
\end{center}
\end{figure}

\begin{figure}[h!]
\centering
\begin{tabular}{cc}
\includegraphics[scale=0.49]{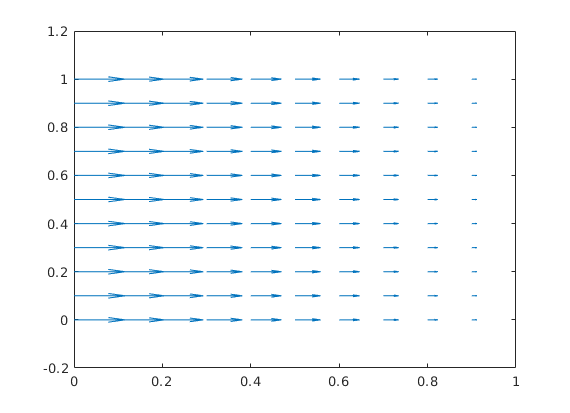} & 
\includegraphics[scale=0.49]{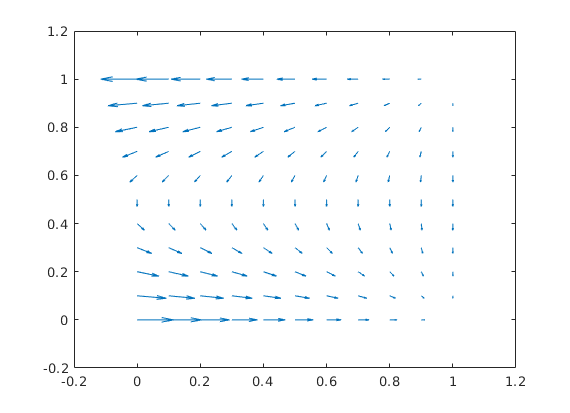} \\
\includegraphics[scale=0.49]{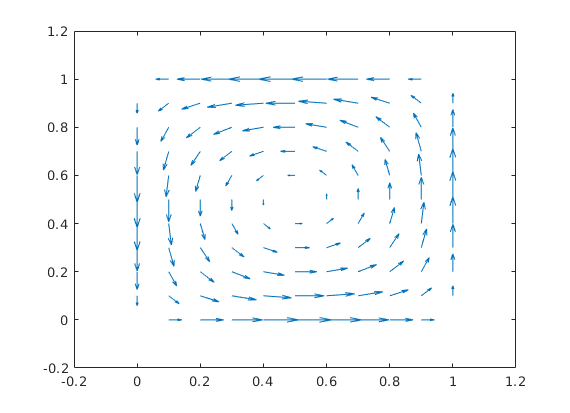} & 
\includegraphics[scale=0.49]{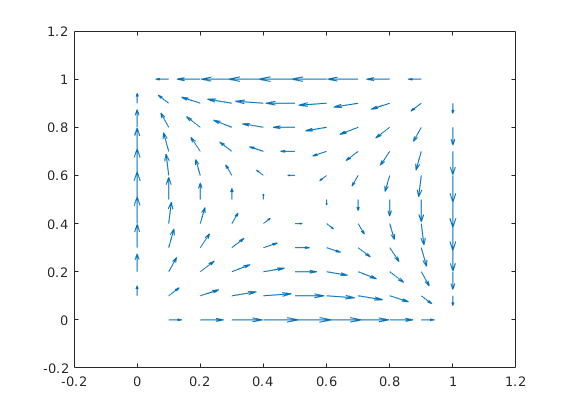} 
\end{tabular}
\caption{Plots of  basis functions on $(0,1)^2$: low order edge function (left above), high order edge based basis function for $p=2$ (right above) to the edge $\mathcal{E}_0$,
high order cell based basis functions for $p=2$ of type 1 and 2 (below).}
\label{fig:Basisfunc}
\end{figure}

With the help of these basis functions, we define the 
two-dimensional N\'ed\'elec space
\begin{equation}
  V_h \coloneqq V^{\mathcal{N}_0}_h\left(\mathcal{T}_h\right)
                \bigoplus_{\mathcal{E} \in \textgoth{E}} V^\mathcal{E}_h\left(\mathcal{T}_h\right)
                \bigoplus_{C \in \mathcal{C}} V^C_h\left(\mathcal{T}_h\right),
\end{equation}
where $V^{\mathcal{N}_0}_h$ is the space of the \textit{lowest-order N\'ed\'elec} function, 
$V^\mathcal{E}_h$ is the space of the \textit{edge-bubbles} and $V^C_h$ is the space of the 
\textit{cell-bubbles}. All basis functions on one element with baryzentric coordinates are displayed in Figure \ref{fig:nedelec_space}.
Visualizations of some basis functions are displayed in Figure \ref{fig:Basisfunc}.
The description of $V_h(\Omega)$ is still not complete, so far we only described $V_h(Q)$,
with $Q$ as previously defined.
It remains to introduce the Piola transformation, which is used to transform the reference
element to any given physical element, see Monk \cite{bk:monk:03} (Lemma 3.57, Corollary 3.58).

\subsection{Discrete weak formulation}
We have gathered everything to write down the discrete weak formulation of the 
time-harmonic Maxwell's equations. We obtain the discrete weak formulation by 
applying the Galerkin method to the equation \eqref{eq:thm_weak}.
Find $E_h \in V_h(\Omega)$ such that
\begin{align}
  \int_\Omega \left( \mu^{-1} \operatorname{curl} \left( \vec E_h \right) \operatorname{curl} \left( \vec \varphi_h \right) 
  - \varepsilon \omega^2 \vec E_h \cdot \vec \varphi_h \right) ~ \dd x 
  + & i \kappa \omega \int_{\Gamma^\infty} \gamma^T \left( \vec E_h \right) \cdot \gamma^T \left( \vec \varphi_h \right) ~ \dd s 
  \nonumber \\ 
  & = \int_{\Gamma^\text{inc}} \gamma^T \left( \vec E^\text{inc} \right) \cdot \gamma^T \left( \vec \varphi_h \right) ~ \dd s
  ~~ \forall \varphi_h \in V_h(\Omega).
\end{align}

\section{Numerical approach}
\label{sec_numerics}
In this section, we first describe domain decomposition and afterwards 
the neural network approximation. In the latter, we also outline 
how to replace the interface operator by the neural network.

\subsection{Domain decomposition}
Since the solution of Maxwell's equation system (\ref{maxwellsystem}) is 
challenging, as already outlined in the introduction, 
we apply a non-overlapping 
domain composition method (DDM)\cite{bk:toselli} in which the domain 
is divided into subdomains as follows
\begin{align*}
  &\overline \Omega=\bigcup_{i=0}^{n_{\mathrm{dom}}} \overline \Omega_i \quad \text{with}\\
  &\Omega_i \cap \Omega_j=\varnothing \quad \forall i \neq j,
\end{align*}
where $n_\mathrm{dom}+1$ is the number of subdomains. 
In such a way, every subdomain $\Omega_i$ becomes small enough so that we 
can handle it with a direct solver.
The global solution of the electric field $E$ is computed via an iterative method, where we solve the time-harmonic Maxwell's
equations on each subdomain with suitable interface conditions between the different subdomains.
Thus, we obtain a solution $E_i^k$ for every subdomain $\Omega_i$, where $k$ denotes the $k$-th iteration step.
The initial interface condition is given by
\begin{align}
  g_{ji}^{k=0}:=-\mu^{-1} \gamma_i^t \left(\mathrm{curl} \left(E_i^{k=0}\right)\right)-i\kappa S\left(\gamma_i^T \left(E_i^{k=0}\right)\right) 
  = 0,
\end{align}
where $S$ describes the interface operator, $i$ is the index of the current domain, and $j$ is the index of the neighbouring domain \cite{art:bour:15}. 
Afterwards, the electric-field $E_i^{k+1}$ is computed at each step by solving the following system
\begin{align} 
  \left\{\begin{array}{lll} 
    \underline{\operatorname{curl}} \left(\mu^{-1} \operatorname{curl} \left( E_i^{k+1} \right) \right) 
      - \omega^2 \varepsilon E_i^{k+1}&=0 & \text{in} \ \Omega_i, \\
    \mu^{-1} \gamma_i^t \left( \operatorname{curl} \left( E_i^{k+1} \right) \right) - i \omega \kappa \gamma_i^T \left( E_i^{k+1} \right) & = 0 & \text{on} \ \Gamma_i^{\infty}, \\ 
    \gamma_i^T\left(E_i^{k+1}\right)&=\gamma_i^T\left(E_i^{\mathrm{inc}}\right) & \text{on} \ \Gamma_i^{\mathrm{inc}}, \\
    \mu^{-1} S\left(\gamma_i^t\left(\operatorname{curl}\left(E_i^{k+1}\right)\right)\right)-i \omega \kappa \gamma_i^T\left(E_i^{k+1} \right)&=g_{ji}^k & \text{on} \ \Sigma_{ij},
  \end{array}\right. 
\label{maxwellsystemDD}
\end{align}
where $\Sigma_{ij}=\Sigma_{ji}:=\partial \Omega_i \cap \partial \Omega_j$ denotes the interface of two neighbouring elements and the interface condition is updated by
\begin{align}
g_{ji}^{k+1}=-\mu^{-1} \gamma_i^t\left(\operatorname{curl}\left(E_i^{k+1}\right)\right)-i\kappa S\left(\gamma_i^T \left(E_i^{k+1}\right)\right)
= -g_{ij}^{k}-2i\kappa S\left(\gamma_i^T \left(E_i^{k+1}\right)\right).
\label{updateinterface}
\end{align}
In case of success we obtain $\lim_{k \to \infty} E_i^k=E|_{\Omega_i}$, but this convergence depends strongly on the chosen interface operator $S$ (see \cite{art:dol09, art:bour:15}). 
The implementation of this approach into deal.II was done 
in \cite{art:beu:21}.

\subsection{Our new approach: Neural network approximation of $S$}
Since the computation of a good approximation of $S$ is challenging, we examine a new approach in which we attempt 
to approximate this operator with the help of a neural network (NN). For a first proof of concept, 
we choose a prototype example and explore whether an NN can approximate the interface values.
As it is not feasible to compute the exact interface operator $S$, we aim to compute $g^{k+l}_{ij},~l>0$ with an NN,
using $g^{k}_{ij}$ and $E^{k+1}_{i}$ as input.
Another benefit of this approach is that we can quickly generate
a training data set from a classical domain decomposition 
method, as described in Section \ref{subsec_training}. 
We choose $S=\mathds{1}$ for simplicity inside our classical domain decomposition method. 
Hence, the advantage of this approach is that one can update the interface condition without recomputing the system (\ref{maxwellsystemDD}) at each step, raising the hope of reducing the computational cost.

\section{Neural network training}
\label{sec_training}
The first step in neural network approximations is the training process,
which is described in this section. Besides the mathematical 
realization, we also need to choose the software libraries. 
We utilize deal.II \cite{software:dealII94} 
to discretize the time-harmonic Maxwell’s equations with the finite element method.
The neural network is trained with PyTorch \cite{software:pytorch}. The exchange of information between the results of the deal.II code and the PyTorch code take place via the hard disk.

\subsection{Basic definitions}
\label{subsec_basics}
First of all, we give a short 
definition of the neural network type employed in this work, and we 
introduce the basic parameters. 
Further information can be found in \cite{Bish06, Cop07, Kri05, BeSha14, Nie15}.
The following notation and descriptions of this subsection are mainly based on \cite{KnoWi21}.

\begin{definition}[Artificial neuron]
An (artificial) neuron (also known as unit \cite{ellacott_1994}, \cite{Bish06}[Section 5.1]) 
  $u$ is a tuple of the form $(\texttt{x}, \texttt{w}, \sigma)$. The components have the following meanings:
\begin{itemize}
\item $\texttt{x} = (\texttt{x}_0,\ldots,\texttt{x}_n)\in \R^{n+1}$ is the input vector. It contains the information, that the neuron receives.
\item $\texttt{w} = (\texttt{w}_0,\ldots,\texttt{w}_n) \in \R^{n+1}$ is the weight vector, 
which determines the influence of the individual input information 
on the output of the neuron. Later, $\texttt{w}$ denotes the weight vector of all neurons.
\item $\sigma : \R \to \R$, with $\sigma = \sum_{i=0}^n \texttt{x}_i \texttt{w}_i \mapsto a$ is the 
activation function. It determines the so-called activation level $a$ from the input and the weights, which represent the output of the neuron.
\end{itemize}
\end{definition}

\begin{definition}[Neural network]
An (artificial) neural (feedforward-) network is a set of neurons
$U$ with a disjoint decomposition $U=U_0 \dot\cup \ldots \dot\cup U_l$. 
The partition sets $U_\mathsf{k}, \mathsf{k}=0,\ldots, l$ are called layers. Here, $U_0$ is the input layer. 
It contains the neurons that receive information from outside. 
Moreover, $U_l$ is the output layer 
with the neurons that return the output. 
Finally, $U_1,\ldots,U_{l-1}$ are the so-called hidden layers. 
\end{definition}
Starting from any neuron $u \in U_\mathsf{k}$, there is a connection to each neuron $\hat{u} \in U_{\mathsf{k}+1}$ for $\mathsf{k}=0,\ldots,l-1$. Such a connection illustrates that the output $a$ of the neuron $u$ is passed on to the neuron $\hat{u}$. This property is the reason for the name feedforward network. 

Each $U_0,\ldots,U_{l-1}$ contains a so-called bias neuron of the form $(0,0,1)$. It has no input, weights and a constant output value $1$ and only transfers a constant bias in the form of the weight to each neuron of the subsequent layer.

\begin{remark}
In the following examples all neurons of the layer $U_\mathsf{k}$ will have the same activation function, given by $\sigma^{(\mathsf{k})}$ for $\mathsf{k}=0,\ldots,l$. Here,
$D_\mathsf{k}:=|U_\mathsf{k}|-1$ for $k=0,\ldots,l-1$ denotes the number of neurons of the $\mathsf{k}$-th layer 
(without bias neuron) and $D_l:=|U_l|$ is the number of neurons of the output layer.
\end{remark}

\subsection{Decomposing the domain}
\begin{figure}[h!]
  \centering
  \begin{tikzpicture}[scale = 5.0]
    \draw (0.5,0.25) node[anchor=south] {\scriptsize $\Omega_0$};
    \draw (0.5,0.75) node[anchor=south] {\scriptsize $\Omega_1$};
    \draw[color=black, thick] (0,0) rectangle (1,1);
    \draw (0.5,0) node[anchor=north] {\scriptsize $\Gamma^{\mathrm{inc}}$};
    \draw (0.5,1) node[anchor=south] {\scriptsize $\Gamma^{\infty}$};
    \draw[color=black, thin] (0,0.5) -- (1,0.5);
    \draw (0.5,0.5) node[anchor=south] {\scriptsize $\Sigma_{01}=\Sigma_{10}$};
    \draw[color=black, thin, ->] (0.15,0.4) -- (0.15,0.6);
    \draw (0.15,0.45) node[anchor=east] {\scriptsize $g_{01}$};
    \draw[color=black, thin, ->] (0.2,0.6) -- (0.2,0.4);
    \draw (0.2,0.55) node[anchor=west] {\scriptsize $g_{10}$};
  \end{tikzpicture}
  \caption{Visualization of the domain $\Omega$ with the chosen decomposition.}
  \label{domain}
\end{figure}
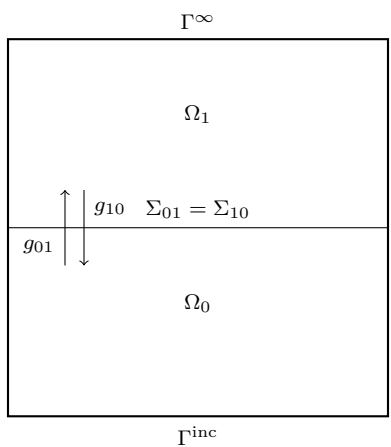

Before constructing the NN, we choose the domain, the decomposition and the grid on 
which the system (\ref{maxwellsystemDD}) is solved to obtain the training values 
because they will influence the network size. The domain in our chosen example, given by
\begin{align*}
  \Omega=(0,1) \times (0,1),
\end{align*}
is divided into two subdomains
\begin{align*}
  \Omega_0=(0,1) \times (0,0.5) \quad \text{and} \quad \Omega_1=(0,1) \times (0.5,1) \quad \text{(see Figure \ref{domain})},
\end{align*}
and the grid on which the FEM 
is applied is a mesh of $32 \times 32$ elements with quadratic N\'ed\'elec elements.
Hence $32$ elements with each $4$ degrees of freedom (dofs) are located on the 
interface in both subdomains. We evaluate the interface condition and the 
solution on each dof and use the values as the NN's input and target. 
Therefore the input contains $4 \cdot \dim(g_{ij})+4 \cdot \dim(E_i)=16$ values, 
and the output consists of $4 \cdot \dim(g_{ji})=8$ values, and we obtain $32$ 
input-target pairs with one computation. 

\subsection{Neural network construction}
\label{subsec_construction}
Regarding the considerations above, we need an input layer with $16$ neurons (without bias) and an 
output layer with $8$ neurons. Furthermore, we use one hidden layer with $500$ neurons (without bias). Hence, for the governing network, we have
\begin{align*}
U=U_0 \dot\cup U_1 \dot\cup U_2
\end{align*}
with $D_0=16$, $D_1=500$ and $D_2=8$. Our tests, presented in Section \ref{sec_tests}, revealed that this is a sufficient size for our purpose. The activation function 
per layer is chosen as follows:
\begin{align*}
\sigma^{(0)}&=id \quad \text{(input layer),} \\
\sigma^{(1)}&=\frac{1}{1+e^{-x}} \quad \text{(hidden layer),} \\ 
\sigma^{(2)}&=id \quad \text{(output layer),}
\end{align*}
where $\sigma^{(1)}$ is known as the sigmoid function, which turned out to be the most effective since the error could be reduced more and more quickly than with other functions we tested e.g.
\begin{align*}
\sigma^{(1)}&=\tanh(x), \\
\sigma^{(1)}&=\log \left(\frac{1}{1+e^{-x}} \right) \quad \text{(LogSigmoid)}, \\
\sigma^{(1)}&=\max(0,x)+\min(0,e^x-1) \quad \text{(CELU)}, \\
\sigma^{(1)}&=a \left( \max(0,x)+\min \left(0,b \left( e^x-1 \right) \right) \right) \quad \text{(SELU)}, \\ 
&\text{with}~a \approx 1.0507,~\text{and}~b \approx 1.6733.
\end{align*}
An exception represents the ReLU function, which we will discuss later in Section \ref{subsec_sigmoid_vs_relu}. Moreover, we apply separate networks $U^{01}$ and $U^{10}$ of the same shape for both interface conditions $g_{01}$
and $g_{10}$ since it turned out that they are approximated differently, fast and accurately. The resulting programming
code is displayed in Figure \ref{fig:code1}.

\begin{figure}[h!]
    \begin{tcolorbox}[breakable, size=fbox, boxrule=1pt, pad at break*=1mm,colback=cellbackground, colframe=cellborder]
\begin{Verbatim}[commandchars=\\\{\}]
\PY{k+kn}{import} \PY{n+nn}{torch}\PY{n+nn}{.}\PY{n+nn}{nn} \PY{k}{as} \PY{n+nn}{nn}
\PY{k+kn}{import} \PY{n+nn}{torch}\PY{n+nn}{.}\PY{n+nn}{nn}\PY{n+nn}{.}\PY{n+nn}{functional} \PY{k}{as} \PY{n+nn}{F}
\PY{k}{class} \PY{n+nc}{Maxwell}\PY{p}{(}\PY{n}{nn}\PY{o}{.}\PY{n}{Module}\PY{p}{)}\PY{p}{:} 
    \PY{k}{def} \PY{n+nf+fm}{\PYZus{}\PYZus{}init\PYZus{}\PYZus{}}\PY{p}{(}\PY{n+nb+bp}{self}\PY{p}{)}\PY{p}{:} 
        \PY{n+nb}{super}\PY{p}{(}\PY{n}{Maxwell}\PY{p}{,} \PY{n+nb+bp}{self}\PY{p}{)}\PY{o}{.}\PY{n+nf+fm}{\PYZus{}\PYZus{}init\PYZus{}\PYZus{}}\PY{p}{(}\PY{p}{)}
        \PY{n+nb+bp}{self}\PY{o}{.}\PY{n}{lin1} \PY{o}{=} \PY{n}{nn}\PY{o}{.}\PY{n}{Linear}\PY{p}{(}\PY{l+m+mi}{16}\PY{p}{,} \PY{l+m+mi}{500}\PY{p}{)} 
        \PY{n+nb+bp}{self}\PY{o}{.}\PY{n}{lin2} \PY{o}{=} \PY{n}{nn}\PY{o}{.}\PY{n}{Linear}\PY{p}{(}\PY{l+m+mi}{500}\PY{p}{,} \PY{l+m+mi}{8}\PY{p}{)} 
        
    \PY{k}{def} \PY{n+nf}{forward}\PY{p}{(}\PY{n+nb+bp}{self}\PY{p}{,} \PY{n}{x}\PY{p}{)}\PY{p}{:} 
        \PY{n}{x} \PY{o}{=} \PY{n}{torch}\PY{o}{.}\PY{n}{sigmoid}\PY{p}{(}\PY{n+nb+bp}{self}\PY{o}{.}\PY{n}{lin1}\PY{p}{(}\PY{n}{x}\PY{p}{)}\PY{p}{)} 
        \PY{c+c1}{\PYZsh{}x = F.relu(self.lin1(x))}
        \PY{n}{x} \PY{o}{=} \PY{n+nb+bp}{self}\PY{o}{.}\PY{n}{lin2}\PY{p}{(}\PY{n}{x}\PY{p}{)} 
        \PY{k}{return} \PY{n}{x}
            
\PY{n}{net} \PY{o}{=} \PY{n}{Maxwell}\PY{p}{(}\PY{p}{)} 
\PY{n+nb}{print}\PY{p}{(}\PY{n}{net}\PY{p}{)}
\end{Verbatim}
\end{tcolorbox}
\caption{PyTorch code of the implementation of the network construction.}
\label{fig:code1}
\end{figure}

\subsection{Training}
\label{subsec_training}
To obtain enough training data, we vary the boundary condition $E^{\mathrm{inc}}$ and create
training and test values to control the network during the training 
and avoid overfitting. The training and test sets are generated by the boundary 
values listed in Table \ref{tab:training}.

\begin{table}[h]
\begin{center}
\begin{tabular}{|c|c|c|}
\hline
\multicolumn{2}{|c|}{$E^{\mathrm{inc}}$ for the training set} & $E^{\mathrm{inc}}$ for the test set \\
\hline
$\left(\begin{array}{c} e^\frac{-(x-0.7)^2}{0.008} \\ 0 \end{array}\right)$ & $\left(\begin{array}{c} \cos(\pi^2 y)+\sin(\pi^2 x) i \\ \sin(\pi^2 y)+0.5 \cos(\pi^2 x) i \end{array}\right)$ & $\left(\begin{array}{c} e^\frac{-(x-0.5)^2}{0.003} \\ 0 \end{array}\right)$ \\ 
\hline
$\left(\begin{array}{c} e^\frac{-(x-0.2)^2}{0.002} \\ 1 \end{array}\right)$ & $\left(\begin{array}{c} \sin(\pi^2 x)+\sin(\pi^2 x) i \\ \sin(\pi^2 y)+0.5 \cos(\pi^2 x) i \end{array}\right)$ & $\left(\begin{array}{c} \cos(\pi^2 y)+\sin(\pi^2 x) i \\ \cos(\pi^2 y)+0.5 \cos(\pi^2 x) i \end{array}\right)$ \\ 
\hline
$\left(\begin{array}{c} e^\frac{-(x-0.7)^2}{0.003} \\ 1 \end{array}\right)$ & $\left(\begin{array}{c} \sin(\pi^2 x)+\sin(\pi^2 x) i \\ \sin(\pi^2 x)+0.5 \cos(\pi^2 x) i \end{array}\right)$ & \\
\hline
$\left(\begin{array}{c} e^\frac{-(x-0.8)^2}{0.003} \\ \sin(\pi^2 x) \end{array}\right)$ & $\left(\begin{array}{c} \cos(\pi^2 y)+\sin(\pi^2 x) i \\ \cos(\pi^2 x)+0.5 \cos(\pi^2 x) i \end{array}\right)$ & \\
\hline
$\left(\begin{array}{c} e^\frac{-(x-0.5)^2}{0.003} \\ \cos(\pi^2 x) \end{array}\right)$ & $\left(\begin{array}{c} \cos(\pi^2 x)+\sin(\pi^2 x) i \\ \cos(\pi^2 y)+0.5 \cos(\pi^2 x) i \end{array}\right)$ & \\
\hline
\end{tabular}
\caption{Boundary values for generating the training set and the test set}
\label{tab:training}
\end{center}
\end{table}

Since we choose $10$ different boundary values for the training set and $2$ for the test set and each of them generates a set of $32$ training/test values 
(one per element on the interface), we obtain all in all a set of $32 \cdot 10=320$ training values and a set of $32 \cdot 2=64$ test values for both networks. To keep the computation simple in a first set of tests, we choose a small wave number $\omega=\frac{2 \pi}{3}$, 
and compute the sets with the iterative DDM in $4$ steps. Afterwards we use the results $\left(g_{ij}^1, E_i^{2}\right)$ as the input and $g_{ji}^3$ as the targets to train our NNs with the application of the mean squared error as the loss function, given by
\begin{align*}
\text{Loss}(\texttt{w})=\frac{1}{2} \sum_{i=1}^N \Vert t^{(i)}-y(\texttt{x}^{(i)},\texttt{w})\Vert^2,
\end{align*}
where $N$ denotes the number of input-target pairs (in our case $N=320$ for the training set and $N=64$ for the test set), $t^{(i)}$ is the target vector, 
$y$ is the function generated by the network and hence $y(\texttt{x}^{(i)},\texttt{w})$ denotes the output of the NN. We refer 
the reader to Section  \ref{subsec_example1} for the specific realization.

As the optimizer, we use the Adam algorithm \cite{kingma2017adam}, which is a line search method based on the following iteration rule
\begin{align*}
x^{\rho+1}=x^\rho+\alpha^\rho p^\rho,
\end{align*}
where $p^\rho$ is called the search direction and $\alpha^\rho$ is the step size (or learning rate in case of NN) for the iteration step $\rho$. The search direction of the Adam algorithm depends on four parameters $\beta_1$, $\beta_2$, $m_1$ and $m_2$, where $\beta_1$ and $\beta_2$ are fixed values in the interval $[0,1)$, and $m_1$ and $m_2$ are updated in each step via
\begin{align*}
m_1^0=m_2^0=0,\quad &m_1^{\rho+1}=\beta_1 m_1^{\rho} + (1-\beta_1) \cdot \nabla \text{Loss}(\texttt{w})^{\rho} \\
\text{and} \quad &m_2^{\rho+1}=\beta_2 m_2^{\rho} + (1-\beta_2) \cdot \Vert \nabla \text{Loss}(\texttt{w})^{\rho} \Vert^2.
\end{align*}
The search direction is then given by
\begin{align*}
p^{\rho-1}=-\widehat{m_1^{\rho}}\Big/\sqrt{\widehat{m_2^{\rho}}+\varepsilon}
\end{align*}
with $\widehat{m_1^{\rho}}=m_1^{\rho}/(1-(\beta_1)^{\rho})$, $\widehat{m_2^{\rho}}=m_2^{\rho}/(1-(\beta_2)^{\rho})$ and $0 < \varepsilon \ll 1$.

The implementation of this training process in PyTorch is displayed in Figure \ref{fig:code2}.

\begin{figure}[H]
\begin{tcolorbox}[breakable, size=fbox, boxrule=1pt, pad at break*=1mm,colback=cellbackground, colframe=cellborder]
\begin{Verbatim}[commandchars=\\\{\}]
\PY{k+kn}{import} \PY{n+nn}{torch}\PY{n+nn}{.}\PY{n+nn}{optim} \PY{k}{as} \PY{n+nn}{optim}
\PY{k+kn}{import} \PY{n+nn}{time}
\PY{n}{start\PYZus{}time} \PY{o}{=} \PY{n}{time}\PY{o}{.}\PY{n}{time}\PY{p}{(}\PY{p}{)}
\PY{n}{tol} \PY{o}{=} \PY{l+m+mf}{3e\PYZhy{}3}
\PY{n}{max\PYZus{}iter} \PY{o}{=} \PY{l+m+mi}{20000}
\PY{n}{iterations} \PY{o}{=} \PY{l+m+mi}{0} 
\PY{n}{loss\PYZus{}test} \PY{o}{=} \PY{n}{tol} \PY{o}{+} \PY{l+m+mi}{1}
\PY{n}{optimizer} \PY{o}{=} \PY{n}{optim}\PY{o}{.}\PY{n}{Adam}\PY{p}{(}\PY{n}{net}\PY{o}{.}\PY{n}{parameters}\PY{p}{(}\PY{p}{)}\PY{p}{,} \PY{n}{lr}\PY{o}{=}\PY{l+m+mf}{1e\PYZhy{}5}\PY{p}{)} 
\PY{n}{criterion} \PY{o}{=} \PY{n}{nn}\PY{o}{.}\PY{n}{MSELoss}\PY{p}{(}\PY{p}{)}
\PY{k}{while}\PY{p}{(}\PY{n}{iterations} \PY{o}{\PYZlt{}} \PY{n}{max\PYZus{}iter} \PY{o+ow}{and} \PY{n}{loss\PYZus{}test} \PY{o}{\PYZgt{}} \PY{n}{tol}\PY{p}{)}\PY{p}{:}
    
    \PY{n}{out} \PY{o}{=} \PY{n}{net}\PY{p}{(}\PY{n}{inp\PYZus{}training}\PY{p}{)}
    \PY{n}{optimizer}\PY{o}{.}\PY{n}{zero\PYZus{}grad}\PY{p}{(}\PY{p}{)}
    \PY{n}{loss} \PY{o}{=} \PY{n}{criterion}\PY{p}{(}\PY{n}{out}\PY{p}{,} \PY{n}{target\PYZus{}training}\PY{p}{)} 
    \PY{n}{loss\PYZus{}test} \PY{o}{=} \PY{n}{criterion}\PY{p}{(}\PY{n}{net}\PY{p}{(}\PY{n}{inp\PYZus{}test}\PY{p}{)}\PY{p}{,} \PY{n}{target\PYZus{}test}\PY{p}{)}
    \PY{n+nb}{print}\PY{p}{(}\PY{l+s+s2}{\PYZdq{}}\PY{l+s+s2}{Loss: }\PY{l+s+si}{\PYZpc{}.5f}\PY{l+s+s2}{\PYZdq{}} \PY{o}{\PYZpc{}} \PY{n}{loss}\PY{p}{,} \PY{l+s+s2}{\PYZdq{}}\PY{l+s+s2}{Test\PYZhy{}Loss: }\PY{l+s+si}{\PYZpc{}.5f}\PY{l+s+s2}{\PYZdq{}} \PY{o}{\PYZpc{}} \PY{n}{loss\PYZus{}test}\PY{p}{)}\PY{p}{)}
    \PY{n}{loss}\PY{o}{.}\PY{n}{backward}\PY{p}{(}\PY{p}{)} 
    \PY{n}{optimizer}\PY{o}{.}\PY{n}{step}\PY{p}{(}\PY{p}{)} 
    \PY{n}{iterations} \PY{o}{+}\PY{o}{=} \PY{l+m+mi}{1}
\PY{n+nb}{print}\PY{p}{(}\PY{l+s+s2}{\PYZdq{}}\PY{l+s+s2}{Final Loss: }\PY{l+s+si}{\PYZpc{}.5f}\PY{l+s+s2}{\PYZdq{}} \PY{o}{\PYZpc{}} \PY{n}{loss\PYZus{}test}\PY{p}{)} 
\PY{n+nb}{print}\PY{p}{(}\PY{l+s+s2}{\PYZdq{}}\PY{l+s+s2}{Number of iterations: }\PY{l+s+si}{\PYZpc{}.0f}\PY{l+s+s2}{\PYZdq{}} \PY{o}{\PYZpc{}} \PY{n}{iterations}\PY{p}{)}
\PY{n}{time\PYZus{}taken} \PY{o}{=} \PY{n}{time}\PY{o}{.}\PY{n}{time}\PY{p}{(}\PY{p}{)} \PY{o}{\PYZhy{}} \PY{n}{start\PYZus{}time}
\PY{n+nb}{print}\PY{p}{(}\PY{l+s+s2}{\PYZdq{}}\PY{l+s+s2}{Run\PYZhy{}Time: }\PY{l+s+si}{\PYZpc{}.4f}\PY{l+s+s2}{ s}\PY{l+s+s2}{\PYZdq{}} \PY{o}{\PYZpc{}} \PY{n}{time\PYZus{}taken}\PY{p}{)}
\end{Verbatim}
\end{tcolorbox}
\caption{PyTorch code of the implementation of the network training.}
\label{fig:code2}
\end{figure}

The network $U^{01}$ is trained with the learning rate $10^{-5}$. The initial training error of $3.12$ and the test error of $5.87$ are reduced to $1.7 \cdot 10^{-4}$ and $3 \cdot 10^{-3}$ after $29\, 843$ training steps. At $U^{10}$, the initial training error of $0.72$ and the test error of $1.28$ are reduced to $3 \cdot 10^{-4}$ and $4 \cdot 10^{-3}$ after $20\, 326$ steps with learning rate of $10^{-5}$ and after further training with a learning rate of $10^{-6}$ in $3706$ steps, we finally achieve the training error $2.9 \cdot 10^{-4}$ and the test error $3 \cdot 10^{-3}$.

\section{Numerical tests}
\label{sec_tests}
In this section, we investigate several numerical experiments
to demonstrate the current capacities of our approach. In addition,
we highlight and analyze shortcomings and challenges.

\subsection{Comparison of new approach and classical DDM}
\label{subsec_example1}
\begin{figure}[H]
\begin{minipage}[b]{.25\linewidth}
\centering
\includegraphics[scale=0.25]{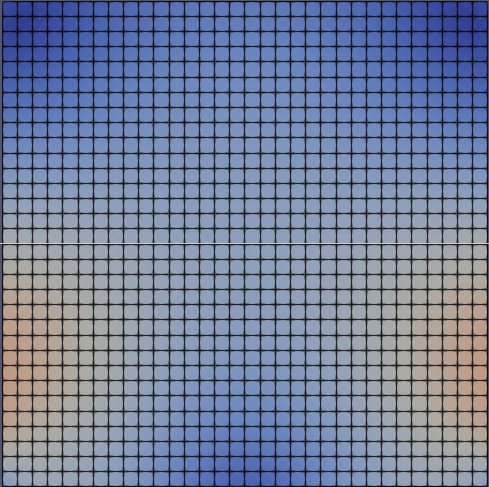}
\end{minipage}
\hspace{.15\linewidth}
\begin{minipage}[b]{.25\linewidth}
\centering
\includegraphics[scale=0.25]{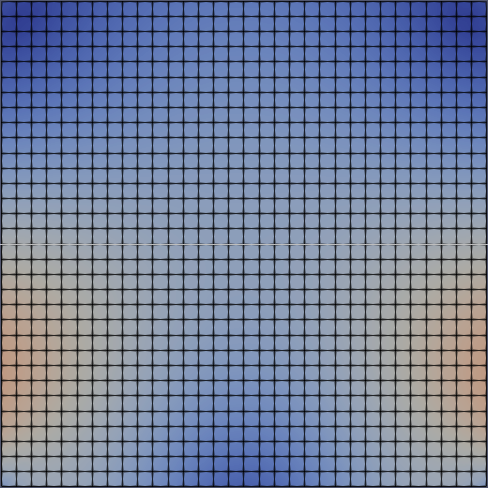}
\end{minipage}
\hspace{.1\linewidth}
\begin{minipage}[b][.36\linewidth][c]{.15\linewidth}
\centering
\includegraphics[scale=0.3]{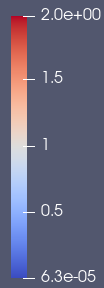}
\end{minipage}
\\\\
\begin{minipage}[b]{.25\linewidth}
\centering
\includegraphics[scale=0.25]{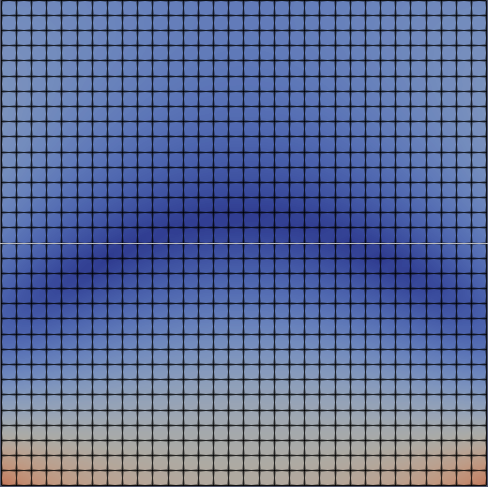}
\end{minipage}
\hspace{.15\linewidth}
\begin{minipage}[b]{.25\linewidth}
\centering
\includegraphics[scale=0.25]{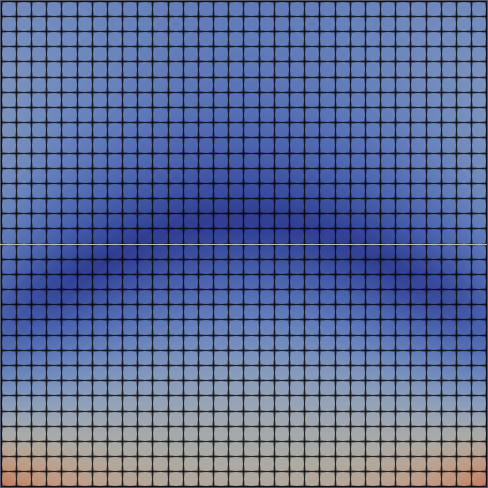}
\end{minipage}
\hspace{.1\linewidth}
\begin{minipage}[b][.36\linewidth][c]{.15\linewidth}
\centering
\includegraphics[scale=0.3]{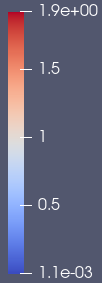}
\end{minipage}
\caption{First example: Real part (above) and imaginary part (below) of the NN solution (left) and the DDM solution (right).}
\label{solution}
\end{figure}

In this first numerical example,
we apply the implemented and trained 
NNs for the following boundary condition
\begin{align*}
  E^{\mathrm{inc}}(x,y)=\left(\begin{array}{c} \cos\left(\pi^2 \left(y-0.5\right)\right)+\sin\left(\pi^2 x\right) i \\ \cos\left(\pi^2 y\right)+0.5 \sin\left(\pi^2 x\right) i \end{array}\right),
\end{align*}
and compute the first interface conditions $g_{10}^1$ and $g_{01}^1$ and the solutions 
$E_1^{2}$ and $E_0^{2}$ by solving (\ref{maxwellsystemDD}) and (\ref{updateinterface}) with the use of the parameters given in Table \ref{tab:parameters}. 
Afterwards, these values are passed on to the networks $U^{01}$ and $U^{10}$. 
The output they return is then handled as our new interface condition, which we use to solve 
system (\ref{maxwellsystemDD}) one more time. With that, we obtain the final solution. 
Moreover, we compute the same example with the DDM in $4$ steps. The results that are 
displayed in Figure \ref{solution} show excellent agreement.

\begin{table}[h]
\centering
\begin{tabular}{|c|c|c|}
\hline
Parameter & Definition & Value \\
\hline
\hline
$\mu$ & relative magnetic permeability & $1.00$ \\
\hline
$\varepsilon$ & relative electric permittivity & $1.49^2$\\
\hline
$\kappa$ & & $\sqrt{\varepsilon}=1.49$ \\
\hline
$\lambda$ & wave length & $3.00$ \\
\hline
$\omega$ & wave number & $\frac{2 \pi}{\lambda}=\frac{2 \pi}{3.00}$ \\
\hline
 & grid size & $\frac{1}{32}$ \\
\hline
\end{tabular}
\caption{Parameters for the DDM}
\label{tab:parameters}
\end{table}

\subsection{Higher wave numbers}
\label{subsec_example2}
\begin{figure}[H]
\begin{minipage}[b]{.25\linewidth}
\centering
\includegraphics[scale=0.25]{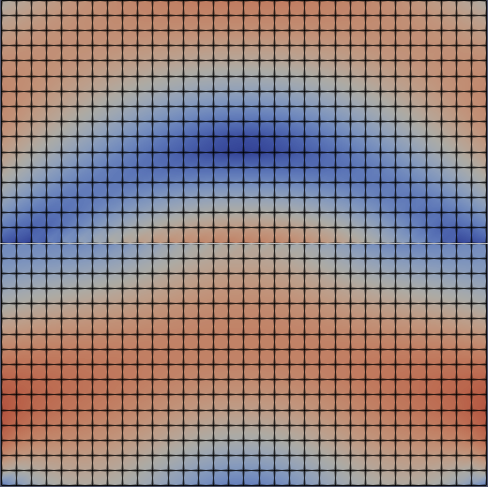}
\end{minipage}
\hspace{.15\linewidth}
\begin{minipage}[b]{.25\linewidth}
\centering
\includegraphics[scale=0.25]{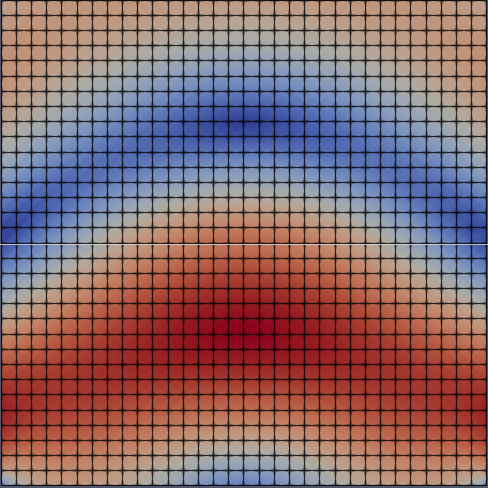}
\end{minipage}
\hspace{.1\linewidth}
\begin{minipage}[b][.36\linewidth][c]{.15\linewidth}
\centering
\includegraphics[scale=0.3]{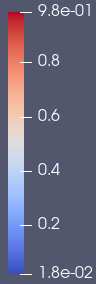}
\end{minipage}
\\\\
\begin{minipage}[b]{.25\linewidth}
\centering
\includegraphics[scale=0.25]{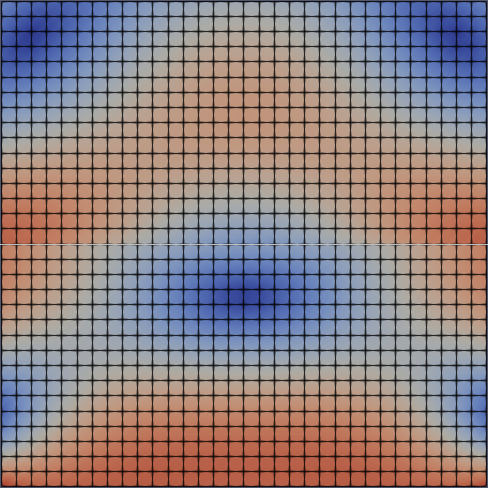}
\end{minipage}
\hspace{.15\linewidth}
\begin{minipage}[b]{.25\linewidth}
\centering
\includegraphics[scale=0.25]{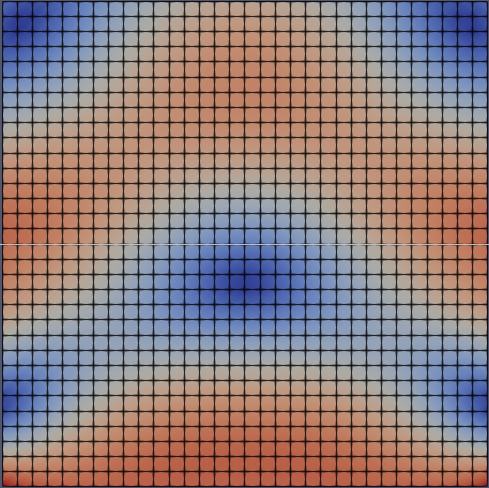}
\end{minipage}
\hspace{.1\linewidth}
\begin{minipage}[b][.36\linewidth][c]{.15\linewidth}
\centering
\includegraphics[scale=0.3]{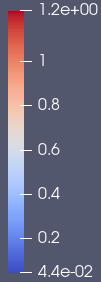}
\end{minipage}
\caption{Second example: Real part (above) and imaginary part (below) of the NN solution (left) and the DDM solution (right).}
\label{solution2}
\end{figure}

As a second example, we increase the wave number, which leads to a more complicated problem.
Therefore we repeat the same computation with $\omega=\pi$ and leave the other parameters 
(especially the parameters and hyperparameters of the neural networks) unchanged. In contrast to the previous example, the results that are displayed 
in Figure \ref{solution2} show differences. While the imaginary part is still well approximated, 
the real part of the NN solution differs significantly from the DDM solution and shows a discontinuity 
on the interface. 

\subsection{Refined computational analysis for intermediate wave numbers}
\begin{figure}[H]
\begin{minipage}[b]{.23\linewidth}
\centering
\includegraphics[scale=0.23]{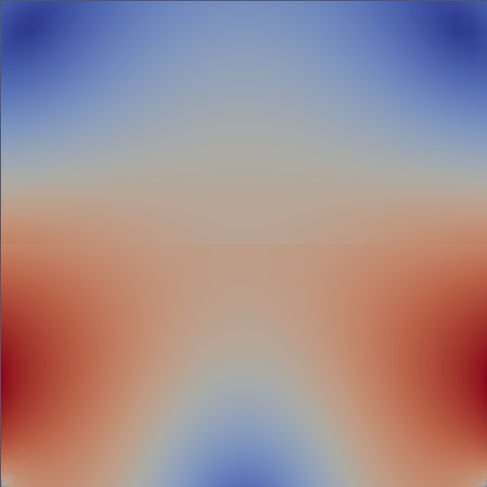}
\end{minipage}
\hspace{.15\linewidth}
\begin{minipage}[b]{.23\linewidth}
\centering
\includegraphics[scale=0.23]{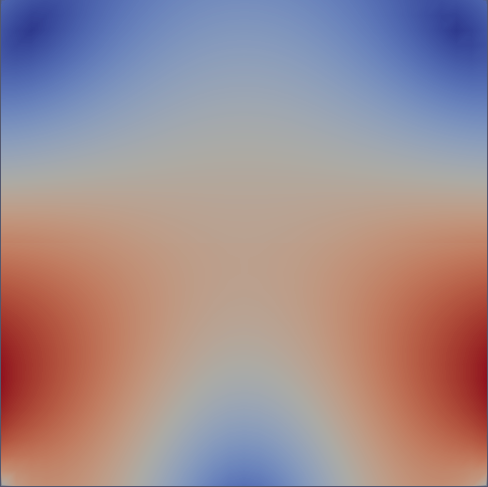}
\end{minipage}
\hspace{.1\linewidth}
\begin{minipage}[b][.25\linewidth][c]{.15\linewidth}
\centering
\includegraphics[scale=0.3]{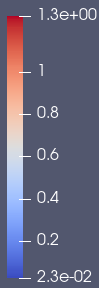}
\end{minipage}
\\\\
\begin{minipage}[b]{.23\linewidth}
\centering
\includegraphics[scale=0.23]{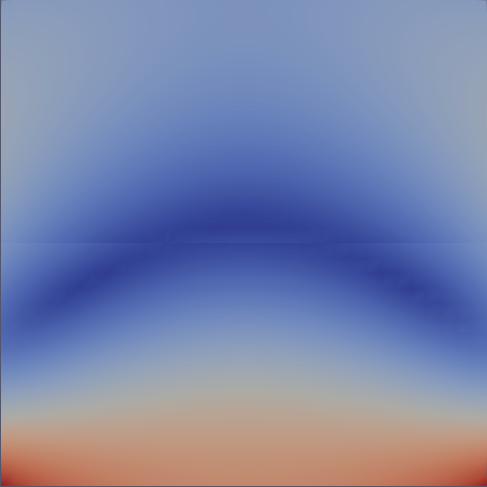}
\end{minipage}
\hspace{.15\linewidth}
\begin{minipage}[b]{.23\linewidth}
\centering
\includegraphics[scale=0.23]{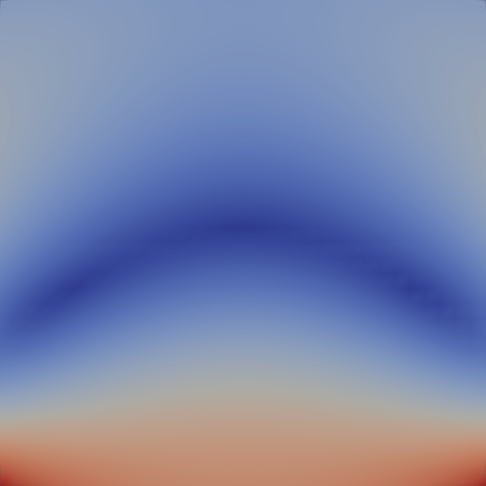}
\end{minipage}
\hspace{.1\linewidth}
\begin{minipage}[b][.25\linewidth][c]{.15\linewidth}
\centering
\includegraphics[scale=0.3]{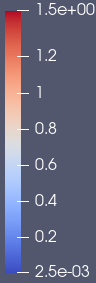}
\end{minipage}
\caption{Third example: Real part (above) and imaginary part (below) of the NN solution (left) and the DDM solution (right) with $\omega=\frac{2 \pi}{2.9}$}
\label{solution3}
\end{figure}

\begin{figure}[H]
\begin{minipage}[b]{.23\linewidth}
\centering
\includegraphics[scale=0.23]{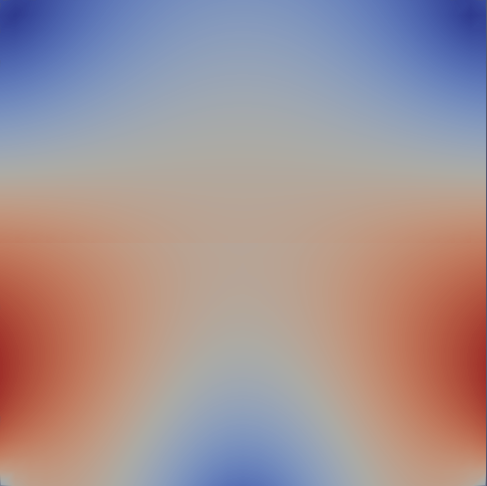}
\end{minipage}
\hspace{.15\linewidth}
\begin{minipage}[b]{.23\linewidth}
\centering
\includegraphics[scale=0.23]{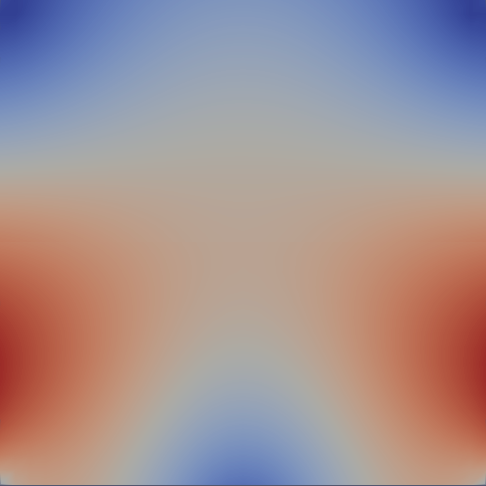}
\end{minipage}
\hspace{.1\linewidth}
\begin{minipage}[b][.25\linewidth][c]{.15\linewidth}
\centering
\includegraphics[scale=0.3]{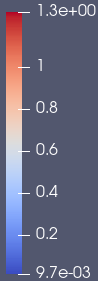}
\end{minipage}
\\\\
\begin{minipage}[b]{.23\linewidth}
\centering
\includegraphics[scale=0.23]{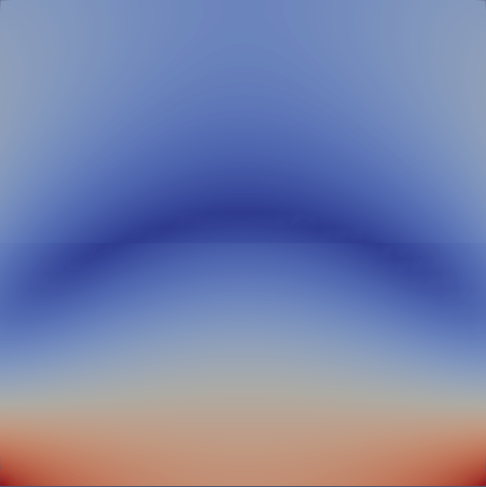}
\end{minipage}
\hspace{.15\linewidth}
\begin{minipage}[b]{.23\linewidth}
\centering
\includegraphics[scale=0.23]{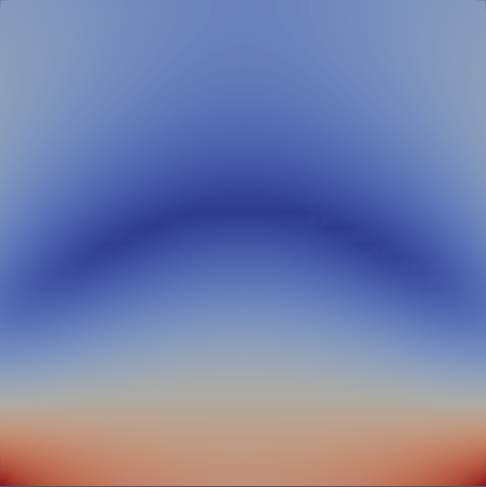}
\end{minipage}
\hspace{.1\linewidth}
\begin{minipage}[b][.25\linewidth][c]{.15\linewidth}
\centering
\includegraphics[scale=0.3]{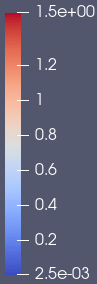}
\end{minipage}
\caption{Fourth example: Real part (above) and imaginary part (below) of the NN solution (left) and the DDM solution (right) with $\omega=\frac{2 \pi}{3.1}$.}
\label{solution4}
\end{figure}

A possible reason for the mismatching results in Section 
\ref{subsec_example2} is the ``problem of big wave numbers'',
which is very well-studied for Helmholtz-type problems \cite{art:gander:12}. The same problem also applies to Maxwell's equations \cite{art:beu:21}.
To verify this conjecture and because of the very distinctive results in 
Section \ref{subsec_example1} and Section \ref{subsec_example2}, 
we attempt two more computations with other wave numbers, namely $\omega=\frac{2 \pi}{2.9}$ and $\omega=\frac{2 \pi}{3.1}$. The results, that are displayed in Figures \ref{solution3} and \ref{solution4}, in which we neglect the representation of the meshes to make the differences more visible, show that the approximation becomes inaccurate if the wave number differs slightly from the one we used for the training, regardless of whether it is larger or smaller. Therefore the bad approximation is not due to the big size of the wave number. Instead of this, it can be assumed that the NNs are specialized for the specific wave number they are trained with and ``learn along'' this value during the training process.

\newpage
\subsection{Comparison of different neural network activation functions: 
Sigmoid vs. ReLU}
\label{subsec_sigmoid_vs_relu}
\begin{figure}[H]
\begin{minipage}[b]{.25\linewidth}
\centering
\includegraphics[scale=0.25]{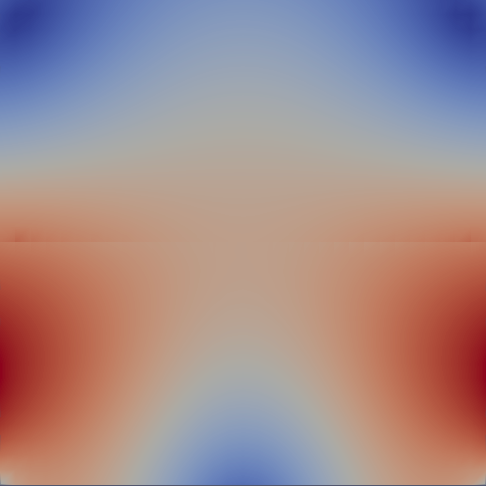}
\end{minipage}
\hspace{.15\linewidth}
\begin{minipage}[b]{.25\linewidth}
\centering
\includegraphics[scale=0.25]{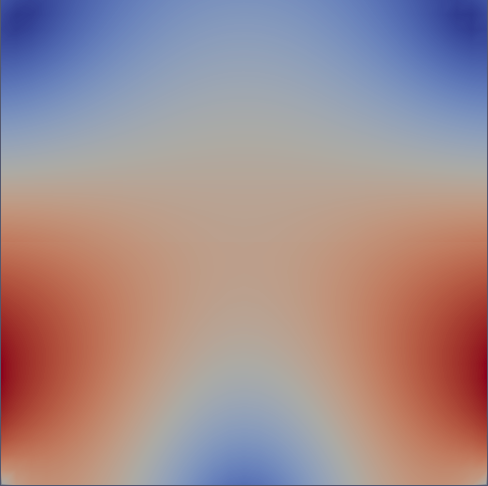}
\end{minipage}
\hspace{.1\linewidth}
\begin{minipage}[b][.36\linewidth][c]{.15\linewidth}
\centering
\includegraphics[scale=0.3]{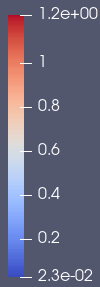}
\end{minipage}
\\\\
\begin{minipage}[b]{.25\linewidth}
\centering
\includegraphics[scale=0.25]{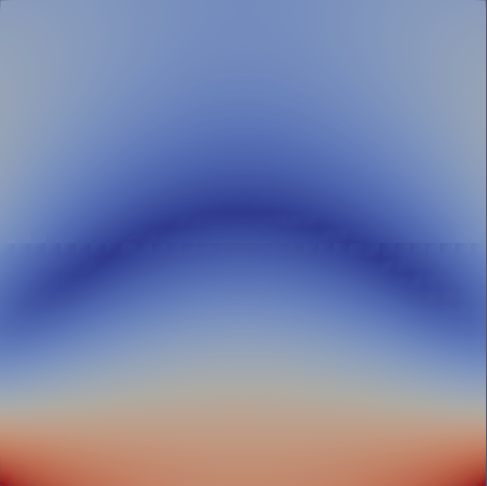}
\end{minipage}
\hspace{.15\linewidth}
\begin{minipage}[b]{.25\linewidth}
\centering
\includegraphics[scale=0.25]{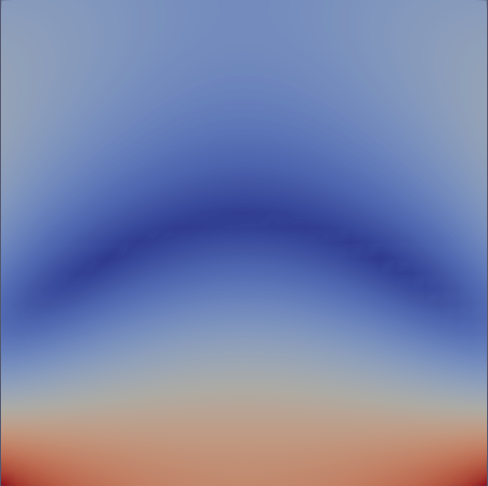}
\end{minipage}
\hspace{.1\linewidth}
\begin{minipage}[b][.36\linewidth][c]{.15\linewidth}
\centering
\includegraphics[scale=0.3]{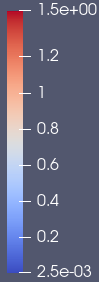}
\end{minipage}
\caption{Real part (above) and imaginary part (below) of the NN solution with the use of ReLU (left) and Sigmoid (right) as activation function.}
\label{fig:sigmoid_vs_relu}
\end{figure}

As mentioned in Section \ref{subsec_construction}, we tested different activation functions to train the NNs before using sigmoid. One of these is the ReLU function given by 
\begin{align*}
f(x)=\max(0,x),
\end{align*}
which is implemented in the PyTorch class torch.nn.functional. This function allows a greater and faster error reduction than the others we tested, including sigmoid. In most cases, the test error of the network $U^{01}$ can be reduced after approx. $16000$ 
steps with a learning rate of $10^{-5}$ and ca. $6500$ 
steps with a learning rate of $10^{-6}$ to $8 \cdot 10^{-4}$, which is almost a quarter compared to the final error in the training of the same NN with sigmoid as the activation function 
(see Section \ref{subsec_training}). 
Also, the test error of $U^{10}$ can be reduced more quickly, 
namely to $2 \cdot 10^{-3}$ after ca. $3500$ 
steps with a learning rate of $10^{-5}$. However, we also observed that the test error grows after a short reduction phase in other cases. But in contrast, the training error continues to shrink, 
revealing that the training of our ReLU-networks is more susceptible to overfitting. This suspicion is strengthened when we apply the successfully trained ReLU-networks to the first example with the same procedure described in Section \ref{subsec_example1}. The results displayed in Figure \ref{fig:sigmoid_vs_relu} show a discontinuity in the interface. This suggests that even in the lucky cases in which the test error is reduced very well, we are dealing with overfitting, and the resulting NNs cannot accurately capture the actual problem. Because of the unreliable training of the ReLU-NNs, it is reasonable to use sigmoid as the activation function instead.

\section{Conclusion}
\label{sec_conclusions}
In this contribution, we provided a proof of concept 
and feasibility study for approximating the interface operator 
in domain decomposition with a feedforward neural network. These 
concepts are applied to the time-harmonic Maxwell's equations. 
We carefully described the numerical framework from the algorithmic 
and implementation point of view. In the realization, we coupled 
deal.II (C++) for solving the Maxwell's equations with PyTorch
for the neural network solution. 
Afterwards, we conducted various numerical tests that included
comparing our new approach with classical domain decomposition.
Then, we studied higher wave numbers in more detail. Therein, we detected 
difficulties, which we further investigated, revealing that the training and testing
of the neural network is highly sensitive to the specific wave number.
Finally, a comparison 
of two different neural network activation functions was undertaken.
As an outlook, we plan to increase the number of subdomains to
study other wave numbers further and apply the method to three-dimensional
Maxwell's equations.

\section*{Acknowledgment}
  This work is funded by the Deutsche Forschungsgemeinschaft (DFG) 
  under Germany’s Excellence Strategy within 
  the Cluster of Excellence PhoenixD (EXC 2122, Project ID 390833453).



\end{document}